# EXISTENCE AND REGULARITY FOR A CLASS OF INFINITE-MEASURE $(\xi, \psi, K)$-SUPERPROCESSES

Alexander Schied

Abstract. We extend the class of $(\xi, \psi, K)$-superprocesses known so far by applying a simple transformation induced by a "weight function" for the one-particle motion. These transformed superprocesses may exist under weak conditions on the branching parameters, and their state space automatically extends to a certain space of possibly infinite Radon measures. It turns out that a number of superprocesses which were so far not included in the general theory fall into this class. For instance, we are able to extend the hyperbolic branching catalyst of Fleischmann and Mueller [FM97] to the case of $\beta$-branching. In the second part of this paper, we discuss regularity properties of our processes. Under the assumption that the one-particle motion is a Hunt process, we show that our superprocesses possess right versions having càdlàg paths with respect to a natural topology on the state space. The proof uses an approximation with branching particle systems on Skorohod space.

## 1. Introduction

In recent years, there has been some effort to construct and describe the most general class of superprocesses. A very general result in this direction is due to Dynkin [Dy94]. He considers the so-called $(\xi, \psi, K)$-superprocesses which are determined by a (time-inhomogeneous) strong Markov process $\xi = (\xi_t, \Pi_{r,x})$ having càdlàg paths in a topological space $(E, \mathcal{B})$, a certain functional $\psi$ acting on the space $\mathcal{B}^+$ of positive measurable functions, and a continuous additive functional $K$ of $\xi$ satisfying

$$(1.1)\qquad \begin{gathered} \sup_x \Pi_{r,x}\big[K[r,t]\big] < \infty, \text{ for all } r < t, \text{ and} \\ \Pi_{r,x}\big[K[r,t]\big] \to 0 \text{ uniformly in } x \text{ as } r,t \to s, \text{ for all } s. \end{gathered}$$

Research supported by Deutsche Forschungsgemeinschaft and, at MSRI, by NSF grant DMS-9701755.

Here and throughout the paper we use the symbols $P$ or $\Pi$ to denote both a probability and the corresponding expectation. The $(\xi,\psi,K)$-superprocess $X=(X_t,P_{r,\mu})$ then is a Markov process taking values in the space $\mathcal{M}$ of finite positive measures on $(E,\mathcal{B})$. It is characterized by the Laplace functionals of its transition probabilities.

$$P_{r,\mu}\big[e^{-\langle f,X_t\rangle}\big]=e^{-\langle v^r,\mu\rangle},\qquad f\in\mathcal{B}^+,$$

where

$$v(r,x)=\Pi_{r,x}\Big[f(\xi_t)-\int_r^t\psi^s(v^s)(\xi_s)\,K(ds)\Big].$$

Even though this class of superprocesses is very large, there are several examples of superprocesses in the literature which are not contained in it. One of the oldest ones is due to Iscoe [Is86], who extended the super-$\alpha$-stable processes to a certain space of *infinite* measures. A similar situation occurs in a recent paper of Dawson and Fleischmann [DF], where $\xi$ is a Brownian motion , $\psi^s(z)(x):=z(x)^2$, and $K$ is the collision local time between $\xi$ and a given typical super-Brownian path started in Lebesgue measure. Condition (1.1) fails in this situation. Also very recently Fleischmann and Mueller [FM97] constructed a superprocess where $\xi$ is one-dimensional Brownian motion killed at 0, $\psi^s(z)(x):=z(x)^2$, and $K(ds)=|\xi_s|^{-\sigma}\,ds$ with $1\le\sigma\le 2$. Again (1.1) fails for $K$, at least if $\sigma=2$.

Our first aim in this paper is to show that these and other superprocesses can be obtained from the class of Dynkin's $(\xi,\psi,K)$-superprocesses by applying a simple transformation. This result can be stated in the general setting of [Dy94], but for simplicity we will restrict ourselves to the case of *local* branching, where $\psi^s(z)(x)$ is a function of $s$, $x$, and $z(x)$. The surprising feature of this transformation is that it allows in many cases to weaken substantially the conditions on the branching characteristics which are imposed in [Dy94]. In the second part of the paper, we focus on the regularity of these superprocesses. It will turn out that the transformed processes have càdlàg sample paths in certain tempered spaces of Radon measures. This is achieved by first showing that a large class of $(\xi,\psi,K)$-superprocesses can be approximated by branching particle systems in the sense of weak convergence on Skorohod space and by proving a lifting property of the quasi-continuous functions of $\xi$ to quasi-continuous linear functions of $X$. Here our assumption that only local branching is allowed becomes crucial. In the recent paper [DFL97] a different approach to weakening (1.1) is taken. By approximation arguments it is shown that, under restrictions on the branching mechanism $\psi$, only the first condition in (1.1) will suffice to construct a $(\xi,\psi,K)$-superprocess. However, also this framework does not include the above examples.

The organization of this paper is the following. In the next section, we state our main results and discuss examples and applications. The existence of our class of processes is shown in Section 3. In Section 4, we collect some facts on branching particle systems, which serve as a preparation for the proof of the weak convergence result. Our regularity theorems are proved in section 6. An appendix on Skorohod spaces closes the paper.

## 2. Statement of the main results

We will mainly use the notation of [Dy94]. Let $\xi = (\xi_t, \mathcal{F}(I), \Pi_{r,x})$ be a (time inhomogeneous) Markov process (cf. [Dy94], Chapter 2) taking values in a metrizable topological space $E$ with Borel field $\mathcal{B}$. We will assume in the sequel that $E$ is a co-Souslin space, which is slightly more general than assuming $E$ to be Polish or Lusinean (cf. [DM75]). Let $b\mathcal{B}$ and $\mathcal{B}^+$ respectively denote the set of all bounded and of all positive $\mathcal{B}$-measurable functions, and let $b\mathcal{B}^+ = \mathcal{B}^+ \cap b\mathcal{B}$. Following [Dy94] the Markov process $\xi$ is called right if $t \mapsto \xi_t(\omega)$ is càdlàg, for all $\omega$, and if $[r,t) \ni s \mapsto \Pi_{s,\xi_s}[f(\xi_t)]$ is $\Pi_{r,\mu}$-a.s. right continuous, for all $f \in \mathcal{B}^+$ and $\mu \in \mathcal{M}$, where $\mathcal{M}$ is the space of all positive finite measures on $(E,\mathcal{B})$ and $\Pi_{r,\mu} := \int \mu(dx)\,\Pi_{r,x}$. Every right process is a strong Markov process (see e.g. [Ku84] or [Dy94]). To avoid measurability problems we will also assume that $(r,x) \mapsto \Pi_{r,x}[f(\xi_t)]$ is measurable, for all $f \in b\mathcal{B}$ and $t > 0$.

If $K$ is a continuous additive functional of $\xi$ and $\psi : \mathbb{R}_+ \times E \times \mathbb{R}_+ \to [0,\infty]$, a $\mathcal{M}$-valued Markov process $X = (X_t, P_{r,\mu})$ is called $(\xi,\psi,K)$-superprocess if

$$P_{r,\mu}\big[e^{-\langle f, X_t\rangle}\big] = e^{-\langle v^r,\mu\rangle}, \qquad f \in \mathcal{B}^+, \tag{2.1}$$

where

$$v(r,x) = \Pi_{r,x}\Big[f(\xi_t) - \int_r^t \psi^s(\xi_s, v^s(\xi_s))\,K(ds)\Big]. \tag{2.2}$$

The following existence result for such superprocesses is taken from [Dy91a], [Dy94]*.

**Theorem 2.1 (Dynkin):** *Suppose $\xi$ is a right process, $K$ is an admissible additive functional of $\xi$ (i.e., $K$ is continuous and satisfies (1.1)), and*

$$\psi^s(x,z) = a^s(x)z + b^s(x)z^2 + \int_{(0,\infty)} (e^{-zu} - 1 + zu)\,\ell^s(x,du), \tag{2.3}$$

*where $a$ and $b$ are positive measurable functions and $\ell$ is a kernel from $\mathbb{R}_+ \times E$ to $(0,\infty)$, and the functions*

$$a^s(x),\ b^s(x),\ \int_{(0,\infty)} u \wedge u^2\,\ell^s(x,du)$$

* [Dy94] assumes that $E$ is a metrizable Luzin space, which is a slightly stronger assumption than ours. However, the crucial results of Section 3.3.4 of [Dy94] also hold under the even weaker assumption that $E$ is Radonian; cf. [Dy91a], [Dy91b].

*are bounded on every set $[0,T]\times E$. Then there exists a $(\xi,\psi,K)$-superprocess, whose transition semigroup is uniquely determined by (2.1) and (2.2), because (2.2) is uniquely solvable if $f\in b\mathcal{B}^+$. Moreover, for $f\in\mathcal{B}^+$,*

$$P_{r,\mu}\big[\langle f,X_t\rangle\big]=\int\int f(y)\,w(r,x,dy)\,\mu(dx),$$

*where $w(r,x,dy)$ is determined by the equations*

$$\int f(y)\,w(r,x,dy)=:w_f(r,x)=\Pi_{r,x}\Big[f(\xi_t)-\int_r^t\int_E w_f(s,\xi_s)\,a^s(\xi_s)\,K(ds)\Big],\qquad f\in b\mathcal{B}^+.$$

Let us now prepare for the statement of our first result. For a $\mathcal{B}$-measurable function $\varrho$ which takes values in $\mathbb{R}_+$, let $\mathcal{B}_\varrho$ denote the set of functions $f$ such that there exists a constant $c$ with $|f|\le c\cdot\varrho$. Dually, $\mathcal{M}^\varrho$ will denote the set of all positive Radon measures $\mu$ on $(E,\mathcal{B})$ with $\langle\varrho,\mu\rangle<\infty$ and $\mu(\varrho=0)=0$. We will call $\varrho$ a *weight function for* $\xi$ if, for every $T>0$, there is a constant $c_T>0$ such that

$$\frac{1}{c_T}\varrho(x)\le\Pi_{r,x}[\varrho(\xi_t)]\le c_T\,\varrho(x),\qquad\text{for all } r\le t\le T.$$

An additive functional $K$ of $\xi$ will be called *$\varrho$-admissible* if it is continuous and satisfies (1.1). Finally, let $\xi^0=(\xi_t,\mathcal{F}(I),\Pi^0_{r,x})$ denote the process obtained by killing $\xi$ when exiting $\{\rho>0\}$.

**Theorem 2.2:** *Suppose $\varrho$ is a weight function for the right process $\xi$, $K$ is $\varrho$-admissible, and $\psi$ is such that $\widetilde{\psi}^s(x,z):=\psi^s(x,z\varrho(x))$ satisfies the assumptions of Theorem 2.1. Then there exists a $\mathcal{M}^\varrho$-valued Markov process $X=(X_t,P_{r,\mu})$ such that*

$$P_{r,\mu}\big[e^{-\langle f,X_t\rangle}\big]=e^{-\langle v^r,\mu\rangle},$$

*where $v$ solves*

$$v(r,x)=\Pi^0_{r,x}\Big[f(\xi_t)-\int_r^t\psi^s(\xi_s,v^s(\xi_s))\,K(ds)\Big].\tag{2.4}$$

*If in addition $f\in\mathcal{B}^+_\varrho:=\mathcal{B}^+\cap\mathcal{B}_\varrho$, then $v\in\mathcal{B}^+_\varrho$ and $v$ is the unique positive solution of (2.4). Moreover, for $f\in\mathcal{B}^+$,*

$$P_{r,\mu}\big[\langle f,X_t\rangle\big]=\int\int f(y)\,w(r,x,dy)\,\mu(dx),$$

*where $w(r,x,dy)$ is determined by the equations*

$$\int f(y)\,w(r,x,dy)=:w_f(r,x)=\Pi^0_{r,x}\Big[f(\xi_t)-\int_r^t w_f(s,\xi_s)\,a^s(\xi_s)\,K(ds)\Big],\qquad f\in\mathcal{B}^+_\varrho,$$

*where $a$ is the linear part in the representation (2.3) corresponding to $\psi$.*

Let us now take a closer look at the applications mentioned in the introduction. We will see that the class of processes constructed in Theorem 2.2 is indeed larger than the one satisfying the assumptions of Theorem 2.1. In particular, we extend the class of hyperbolic superprocesses of [FM97] to the case of $\beta$-branching. By following the recipes of the examples below one can easily get many new superprocesses.

**Examples: (i) Hyperbolic superprocesses:** Let $\xi$ be a one-dimensional Brownian motion stopped at the first hitting time $\tau$ of 0. Then $\varrho(x) := |x|$ is harmonic for $\xi$. In particular, $\varrho$ is a weight function. Let $0 < \beta \le 1$, $1 \le \sigma \le 1+\beta$, and

$$K(ds) := 1 \vee |\xi_s|^{1+\beta-\sigma} \mathrm{I}_{\{s\,<\,\tau\}}\, ds$$

Using $K(ds) \le (\mathrm{I}_{\{s\,<\,\tau\}} + \rho(\xi_s))\, ds$ and the reflection principle for Brownian motion we easily get that

$$\begin{aligned} \Pi_{r,x}\big[K[r,t]\big] &\le \Pi_{r,x}\big[(\tau - r) \wedge (t-r)\big] + \rho(x)(t-r) \\ &\le \rho(x)\Big(\sqrt{\frac{2(t-r)}{\pi}} + t - r\Big). \end{aligned}$$

Thus $K$ is $\varrho$-admissible. The function

$$\psi^s(x,z) := \Big(\frac{z}{|x|^{\sigma/(1+\beta)} \vee |x|}\Big)^{1+\beta}$$

satisfies the assumptions of Theorem 2.2. Hence there exists an $\mathcal{M}^\rho$-valued $(\xi^0, \psi, K)$-superprocess associated with the equation

$$v(r,x) = \Pi^0_{r,x}\Big[f(\xi_t) - \int_r^t v(s,\xi_s)^{1+\beta} \frac{1}{|\xi_s|^\sigma}\, ds\Big], \qquad f \in \mathcal{B}^+_{\tilde\rho},$$

if $0 < \beta \le 1$ and $1 \le \sigma \le 1+\beta$. For $\beta = 1$, this gives the class of processes constructed in [FM97] (note that there only *finite* starting measures were considered). It is worth remarking that in all these models one has $P_{r,\delta_x}\big[\langle 1, X_t\rangle\big] = \Pi^0_{r,x}[\xi_t \ne 0] \to 0$ as $t \uparrow \infty$. This was also observed in Theorem 3 (b) of [FM97].

**(ii) Infinite-measure super-$\alpha$-stable processes:** Let $\xi$ be a standard Brownian motion on $E = \mathbb{R}^d$. Then, if $p > d$, $\phi_p(x) := (1+|x|^2)^{-p/2}$ is a weight function for $\xi$. If $\xi$ is $\alpha$-stable with $\alpha \in (0,2)$, then $\phi_p$ is a weight function if $p < d + \alpha$. See [Is86], Corollary 2.4 and its proof for both results. Now let $K(ds) = \phi_p(\xi_s)^{1+\beta}\, ds$, for $0 < \beta \le 1$. Then

$$\phi_p(x)^{-1}\Pi_{r,x}\big[K[r,t]\big] \le \int_r^t \phi_p(x)^{-1}\Pi_{r,x}[\phi_p(\xi_s)]\, ds \le c_t(t-r).$$

Thus $K$ is $\phi_p$-admissible. Next let $\psi^s(x,z) = \big(z/\phi_p(x)\big)^{1+\beta}$. Then, if $X = (X_t, P_{r,\mu})$ is the $\mathcal{M}^{\phi_p}$-valued $(\xi,\psi,K)$-superprocess of Theorem 2.2, $v(r,x) := -\log P_{r,\delta_x}\big[e^{-\langle f, X_t\rangle}\big]$ solves

$$v(r,x) = \Pi_{r,x}\Big[f(\xi_t) - \int_r^t v(s,\xi_s)^{1+\beta}\,ds\Big].$$

Thus we have found another construction of Iscoe's infinite-measure superprocesses.

**(iii) Superprocesses with locally admissible branching functional:** Let $\xi$ be Brownian motion on $\mathbb{R}^d$ and $\phi_p$ as in (i) with $p > d$. Then let $K$ be such that

$$\widetilde{K}(ds) := \phi_p^2(\xi_s)\,K(ds) \text{ is } \phi_p\text{-admissible.} \tag{2.5}$$

By taking $\psi^s(x,z) = \big(z/\phi_p(x)\big)^2$ we see that our $(\xi,\psi,\widetilde{K})$-superprocess is associated with the equation

$$v(r,x) = \Pi_{r,x}\Big[f(\xi_t) - \int_r^t v(s,\xi_s)^2\,K(ds)\Big]. \tag{2.6}$$

Thus our Theorem 2.2 yields a variant of Proposition 12 (a) of [DF97]. There existence of a superprocess associated with (2.6) has been shown under the slightly different condition that $\phi_p(\xi_s)\,K(ds)$ is admissible for $\xi$ in the sense of (1.1). This property is called *local admissibility* in [DF97]. Together with (2.5) local admissibility implies existence of second moments of the $(\xi,\psi,\widetilde{K})$-superprocess (Proposition 12 (b) of [DF97]). Of course, our construction immediately extends to the case of $\beta$-branching with the appropriate modifications.

Next we are interested in the regularity of the process constructed in Theorem 2.2. To this end, recall first that a right Markov process $\xi$ is called *Hunt process* if it quasi-left continuous. I.e., $\xi_{\sigma_n} \to \xi_\sigma$ $\Pi_{r,x}$-a.s., for every increasing sequence $(\sigma_n)$ of $(\mathcal{F}[r,t])_{t\geq r}$-stopping times such that $\sigma = \lim_n \sigma_n$.

Now endow $\mathcal{M}^\varrho$ with the topology generated by the mappings

$$\mu \mapsto \langle f,\mu\rangle, \qquad \text{for continuous } f \in \mathcal{B}_\varrho.$$

We will call this the *$\varrho$-weak topology* on $\mathcal{M}^\varrho$. The case $\varrho \equiv 1$ then gives the usual weak topology on $\mathcal{M}$. Note that we do not assume that $\varrho$ itself is continuous.

**Theorem 2.3:** *In addition to the assumptions of Theorem 2.2 suppose that $\xi$ is a Hunt process. Then $X_t$ has a càdlàg version with respect to the $\varrho$-weak topology, and $X$ is a right process.*

In the case $\varrho \equiv 1$ and $K(ds) = ds$, Fitzsimmons [Fi88] proved that $X$ is a Hunt process with respect to the Ray weak* topology on $\mathcal{M}$. Kuznetsov [Ku94] proved the strong Markov property and the existence of a right continuous version for the class of superprocesses described in Section 5.5.3 of [Dy94]. Dawson, Fleischmann and Leduc [DFL97] extended Fitzsimmon's result to superprocesses with branching functionals $K$ having bounded characteristic (i.e., the first condition in (1.1) holds), under the additional assumptions that $E$ is compact, $\xi$ is a Feller process, and $\psi$ corresponds to finite variance branching.

The existence of a *càdlàg* version in the generality of Theorem 2.3 seems to be new even for $\rho \equiv 1$. In the latter case, we will prove it by showing that the $(\xi, \psi, K)$-superprocess can be approximated by branching particle systems with respect to weak convergence on Skorohod space. This relies on a lifting of the quasi-continuous functions $f(r, x)$ of $\xi$ to quasi-continuous linear functions $\langle f^r, \mu\rangle$ of $X$. This lifting is then also applied to the transformation we use in the case of a general weight function $\varrho$.

A *branching particle system* is a Markov process $X = (X_t, P^1_{r,\nu})$ taking values in the space $\mathcal{M}_0 \subset \mathcal{M}$ of integer valued measures. Each particle moves independently from all others according to the right Markov process $\xi$, and dies after a random time. Given the path of the particle, the probability of surviving the time interval $[r, t]$ is given by $e^{-K[r,t]}$, where $K$ is an admissible additive functional of $\xi$. At the time of its death, the particle produces a random number of children according to a probability $q(t, x, dn)$ on $\{0, 1, \ldots\}$. Here $t$ and $x$ stand for the time and spatial location of the death of the particle respectively. Afterwards the offspring perform the same process, starting independently from $(t, x)$. It is well-known that (at least under condition (2.11) below)

$$P^1_{r,\nu}\Big[e^{-\langle f, X_t\rangle}\Big] = e^{\langle \log u^r, \mu\rangle}, \qquad f \in \mathcal{B}^+,\ \nu \in \mathcal{M}_0, \tag{2.7}$$

where $u$ solves

$$u(r, x) = \Pi_{r,x}\Big[e^{-f(\xi_t)} + \int_r^t \widetilde{\varphi}^s(\xi_s, u^s(\xi_s))\, K(ds)\Big], \tag{2.8}$$

with

$$\widetilde{\varphi}(s, x, z) := \varphi(s, x, z) - z := \int z^n\, q(s, x, dn) - z.$$

See [Dy91a], [Dy94]. Next suppose that, for every $\beta \in (0, 1]$, we are given a kernel $q_\beta(s, x, dn)$ or equivalently the corresponding generating function $\varphi_\beta$, and replace $K$ by $\frac{1}{\beta}K$. The corresponding branching particle system will be denoted by $P^\beta_{r,\nu}$. For $\mu \in \mathcal{M}$ arbitrary, let $\pi_{\mu/\beta}$ denote the Poisson measure on $\mathcal{M}_0$ with intensity $\mu/\beta$. Then

$$P^\beta_{r,\pi_{\mu/\beta}}\Big[e^{-\langle \beta f, X_t\rangle}\Big] = e^{-\langle v^r_\beta, \mu\rangle}, \qquad f \in \mathcal{B}^+, \tag{2.9}$$

where

$$v_\beta(r,x) = \Pi_{r,x}\Big[\frac{1-e^{-\beta f(\xi_t)}}{\beta} + \int_r^t \psi^s_\beta(\xi_s, v^s_\beta(\xi_s))\,K(ds)\Big], \tag{2.10}$$

and

$$\psi^s_\beta(x,z) = \frac{1}{\beta^2}\Big(\varphi^s_\beta(x,1-\beta z) - 1 + \beta z\Big).$$

Our next result is taken from [Sc92]. It generalizes previous work by Watanabe [Wa68], Ethier and Kurtz [EK86], Roelly [R-C86], and Dawson and Perkins [DP91], and it complements the fdd-convergence proved in [Dy94], [Dy91a].

**Theorem 2.4:** *Suppose that $\xi$ is a Hunt process and that $K$ is admissible for $\xi$. Assume in addition that*

$$\sup_{s,x} \int n\, q_\beta(s,x,dn) \le 1 \qquad \textit{for all } \beta \in (0,1], \tag{2.11}$$

*and that there exists a function $\psi$ satisfying the assumptions of Theorem 2.1 such that $\psi^s_\beta(x,z) \to \psi^s(x,z)$ uniformly on each set $[0,T] \times E$. Then the distributions of $\beta X_t$ under $P^\beta_{r,\pi_{\mu/\beta}}$ converge weakly on the Skorohod space $D(\mathbb{R}_+,\mathcal{M})$ to the $(\xi,\psi,K)$-superprocess started in $(r,\mu)$.*

For each $\psi$ satisfying the assumptions of Theorem 2.1, [Dy91a] gives an approximation of branching particle systems such that the assumptions of Theorem 2.4 hold.

For the proof of Theorem 2.4, we will need the also otherwise useful *historical process*, which was introduced in [DP91]. The *path process* associated with $\xi$ is the $\widehat{E} := D(\mathbb{R}_+;E)$-valued process $\widehat{\xi}_t$ defined by $\widehat{\xi}_t(s) = \xi_{t\wedge s}$. It is then easy to construct a Markov process $\widehat{\xi} = (\widehat{\xi}_t, \mathcal{F}[r,t], \widehat{P}_{r,\hat{x}})$ such that $\widehat{\xi}_r = \widehat{x}(\cdot\wedge r)$ $\widehat{P}_{r,\hat{x}}$-a.s. and such that the law of $[r,\infty) \ni s \mapsto \widehat{\xi}_t(s)$ under $\widehat{P}_{r,\hat{x}}$ equals the law of $[r,\infty) \ni s \mapsto \xi_{t\wedge s}$ under $P_{r,\hat{x}(r)}$; see e.g. [Dy91b].

**Remark:** (i) $\widehat{E} = D(\mathbb{R}_+;E)$ is again a metrizable co-Souslin space (see Proposition 7.1 (ii) below). This is in fact our reason for working with this class of phase spaces instead with the more common class of metrizable Lusin spaces. Indeed, Kolesnikov [Ko99] shows that $D(\mathbb{R}_+;\mathbb{Q})$ is *not* a Lusin space even though $\mathbb{Q}$, the set of rationals, is.

(ii) Unless $\widehat{x} = \widehat{x}(\cdot\wedge r)$ the path process $\widehat{\xi}$ lacks the normal property $\widehat{P}_{r,\hat{x}}[\widehat{\xi}_r = \widehat{x}] = 1$, which is assumed in [Dy94]. However, one can check that it is not needed for the proof of Theorem 2.1. Alternatively, one can work with the space-time process $(t,\widehat{\xi}_t)$ and project down afterwards.

Let $\widehat{\psi}^s(\widehat{x},z) := \psi^s(\widehat{x}(s),z)$ with $\psi$ as in Theorem 2.1. The $(\widehat{\xi},\widehat{\psi},K)$-superprocess is usually called historical superprocess associated with the $(\xi,\psi,K)$-superprocess. Analogously one can define historical branching particle systems. If $\pi_t : D(\mathbb{R}_+;E) \to E$ is given by $\pi_t(\widehat{x}) = \widehat{x}(t)$ one can recover (the finite dimensional marginals of) the $(\xi,\psi,K)$-superprocess $X$ from the $(\widehat{\xi},\widehat{\psi},K)$-superprocess $\widehat{X}$ by setting $X_t = \widehat{X}_t \circ \pi_t^{-1}$. Even though $\pi_t$ is not continuous, we have the following result.

**Corollary 2.5:** *(i) Suppose that the $(\xi,\psi,K)$-superprocess is approximated by branching particle systems as in Theorem 2.4. Then the corresponding historical branching particle systems converge weakly on path space to the corresponding historical superprocess.*

*(ii) If $\widehat{X}_t$ is a càdlàg version of the $(\widehat{\xi},\widehat{\psi},K)$-superprocess, then $t \mapsto X_t = \widehat{X}_t \circ \pi_t^{-1}$ is a càdlàg version of the $(\xi,\psi,K)$-superprocess.*

**Remark:** In some cases, it is possible to apply (a variant of) Theorem 2.2 also for historical superprocesses. For instance, if $E = \mathbb{R}^d$, $\varrho = \phi_p$, and $\xi$ is Brownian motion as in the example above, then Lemma 7 of [Sc96] implies that

$$\phi_p(x) \le \Pi_{r,x}\Big[\sup_{s\le t} \phi_p(\xi_s)\Big] \le c\phi_p(x).$$

So, with the appropriate modifications, one can use $\widehat{\phi}_p(\widehat{x}) := \sup_t \phi_p(\widehat{x}(t))$ as a weight function.

## 3. Proof of Theorem 2.2

First we show existence up to a fixed given time $T > 0$. To this end, let $h(r,x) = \Pi_{r,x}[\varrho(\xi_T)]$. By assumption it follows that

$$(3.1)\qquad \varrho(x) = 0 \iff h(t,x) = 0\ \forall t \le T \iff \exists\, t \le T \text{ such that } h(t,x) = 0.$$

The process $t \mapsto h(t,\xi_t)$ is a positive right continuous martingale up to $T$. In particular,

$$(3.2)\qquad h(t,\xi_t) = 0, \text{ for all } t \in [r,T],\ \Pi_{r,x}\text{-a.s. if } h(r,x) = 0.$$

Let $\xi^h = (\xi_t, \Pi^h_{r,x})$ denote the corresponding $h$-process. I.e.,

$$\Pi^h_{r,x}[A] = \frac{1}{h(r,x)}\Pi_{r,x}\big[\varrho(\xi_T);\, A\big], \qquad A \in \mathcal{F}[r,T],$$

if $h(r,x) > 0$, and $\Pi^h_{r,x} := \Pi_{r,x}$ otherwise. It is then easy to check that $\xi^h$ is again a right Markov process.

Note the the $\varrho$-admissibility of $K$ implies that $K[r,t]=0$, for all $t \geq r$, $\Pi_{r,x}$-a.s. if $h(r,x)=0$. Hence the strong Markov property of $\xi$ implies that

$$K[\tau,\infty)=0 \qquad \Pi_{r,x}\text{-a.s., for all } r,x, \tag{3.3}$$

if $\tau=\inf\{t\,|\,\varrho(\xi_t)=0\}$. Therefore we can define a continuous additive functional $K^h$ of $\xi^h$ by

$$K^h(ds)=\frac{1}{h(s,\xi_s)}\,K(ds).$$

Then the right property of $\xi^h$ and the continuity of $K$ imply that

$$\Pi^h_{r,x}\big[K^h[r,t]\big]=\frac{\Pi_{r,x}\big[K[r,t]\big]}{h(r,x)}, \tag{3.4}$$

as can be seen by approximating $K[r,t]=\int_r^t h(s,\xi_s)\,K^h(ds)$ by Riemann sums and by applying the Markov property of $\xi$. Hence it follows that $K^h$ is admissible for $\xi^h$. Next let $\psi^s_h(x,z):=\psi^s(x,h^s(x)z)$. One easily checks that $\psi_h$ satisfies the assumptions of Theorem 2.1. Thus there exists a $(\xi^h,\psi_h,K^h)$-superprocess $X^h=(X_t,P^h_{r,\mu})$.

Now we claim that

$$X_t(\varrho=0)=0,\ P^h_{r,\mu}\text{-a.s., for all } r\leq t \text{ and } \mu\in\mathcal{M} \text{ with } \mu(\varrho=0)=0. \tag{3.5}$$

Indeed, $P^h_{r,\delta_x}[X_t(\varrho=0)]=:w(r,x)$ is determined by

$$w(r,x)=\Pi^h_{r,x}\big[\varrho(\xi_t)=0\big]-\Pi^h_{r,x}\Big[\int_r^t w(s,\xi_s)h^s(\xi_s)a^s(\xi_s)\,K^h(ds)\Big]. \tag{3.6}$$

Note that it follows from (3.1), (3.2), and our definition of $\Pi^h_{r,x}$ that

$$\Pi^h_{r,x}\big[\varrho(\xi_t)=0\big]=\mathrm{I}_{\{h(r,x)=0\}}.$$

Hence (3.3) implies that $w(r,x):=\mathrm{I}_{\{h(r,x)=0\}}$ solves (3.6). Therefore (3.5) follows.

Since, for $f\in\mathcal{B}^+$, $v_h(r,x):=-\log P^h_{r,\delta_x}\big[e^{-\langle f/h^t,X_t\rangle}\big]$ solves

$$\begin{aligned} v_h(r,x)&=\Pi^h_{r,x}\Big[\frac{f(\xi_t)}{h(t,\xi_t)}-\int_r^t\psi^s_h(\xi_s,v^s_h(\xi_s))\,K^h(ds)\Big]\\ &=\frac{1}{h(r,x)}\Pi^0_{r,x}\Big[f(\xi_t)-\int_r^t\psi^s_h(\xi_s,v^s_h(\xi_s))\cdot h^s(\xi_s)\,K^h(ds)\Big], \end{aligned} \tag{3.7}$$

it follows that $v(r,x) := h(r,x)\cdot v_h(r,x)$ solves

$$v(r,x) = \Pi^0_{r,x}\Big[f(\xi_t) - \int_r^t \psi^s(\xi_s, v^s(\xi_s))\, K(ds)\Big].$$

Moreover, if $f \in \mathcal{B}^+_\varrho$, then $f/h^t \in b\mathcal{B}^+$ and thus $v_h \in b\mathcal{B}^+$ which in turn implies $v \in \mathcal{B}^+_\varrho$.

Now suppose that $\widetilde{v}$ is another positive solution of (2.4), with $f \in \mathcal{B}_\varrho$. Then $f/h^t \in b\mathcal{B}^+$. But also $\widetilde{v}$ is of the form $\widetilde{v}(r,x) = v'(r,x)\cdot h(r,x)$ with $v' \in \mathcal{B}^+$, because (3.1) and (3.3) imply that any solution of (2.4) vanishes on $\{(r,x)\,|\,h(r,x)=0\}$. Going back the above steps one finds that $v'$ solves (3.7), and hence $v' = v_h$ by uniqueness for (3.7). Thus the positive solution of (2.4) is unique if $f \in \mathcal{B}^+_\varrho$.

Next suppose we are given $\mu \in \mathcal{M}^\varrho$. Define an element $\mu^h$ of $\mathcal{M}$ by $\mu^h(dx) = h(r,x)\,\mu(dx)$. Due to (3.5) we can define $P_{r,\mu}$ to be the distribution of $h(t,x)^{-1}X_t(dx)$ under $P^h_{r,\mu^h}$ up to time $T$. By construction this process lives in $\mathcal{M}^\varrho$. Moreover, the facts that, for $f \in pB_\varrho$, $v$ is unique and $v^r \in \mathcal{B}^+_\varrho$, imply the Markov property and the uniqueness of the finite dimensional distributions of $P_{r,\mu}$. In particular, our construction is independent of $T$, and we can extend our process from $[r,T]$ to $[r,\infty)$.

To prove the first moment formula note that by construction and (3.5)

$$P_{r,\mu}\Big[\langle f, X_t\rangle\Big] = P_{r,\mu}\Big[\langle f\cdot \mathrm{I}_{\{\varrho>0\}}, X_t\rangle\Big] = \int h(r,x)\cdot w^h_f(r,x)\,\mu(dx),$$

where $w^h_f(r,x)$ solves

$$\begin{aligned} w^h_f(r,x) &= \Pi^h_{r,x}\Big[\frac{f(\xi_t)}{h(t,\xi_t)}\Big] - \Pi^h_{r,x}\Big[\int_r^t w^h_f(s,\xi_s)h^s(\xi_s)a^s(\xi_s)\,K^h(ds)\Big] \\ &= \frac{1}{h(r,x)}\Big(\Pi^0_{r,x}\big[f(\xi_t)\big] - \Pi^0_{r,x}\Big[\int_r^t w^h_f(s,\xi_s)h^s(\xi_s)a^s(\xi_s)\,K(ds)\Big]\Big), \end{aligned}$$

where we have used (3.2). Thus $w_f(r,x) := P_{r,\delta_x}\big[\langle f, X_t\rangle\big]$ solves

$$w_f(r,x) = \Pi^0_{r,x}\big[f(\xi_t)\big] - \Pi^0_{r,x}\Big[\int_r^t w_f(s,\xi_s)a^s(\xi_s)\,K(ds)\Big]$$

We leave it to the reader to check that $w(r,x,dy)$ is in fact uniquely determined by these equations. $\square$

## 4. Preliminary remarks on branching particle systems

With $(X_t, P^1_{r,\nu})$ we will denote the stochastic process modeling our branching particle system with generating function $\varphi$ subject to condition (2.11). Most results in this section are standard, so we only sketch the proofs.

**Lemma 4.1:** *Suppose $K$ is admissible. Then, for any $T > 0$ and every $\varepsilon > 0$, there is a $\delta > 0$ such that*

$$k(r,t) := \sup_{r\le s\le t}\sup_x \Pi_{s,x}\big[K[s,t]\big] < \varepsilon, \tag{4.1}$$

*for all $0 \le r < t \le T$ with $t - r < \delta$.*

**Proof:** By Lemma 3.3.1 of [Dy94] there is a partition $0 = t_0 < t_1 < \cdots < t_n = T$ such that $k(t_{i-1}, t_i) < \varepsilon/2$, for $i = 1, \dots, n$. The assertion then follows by taking $\delta$ as $\min\{t_i - t_{i-1} \,|\, i = 1, \dots, n\}$. □

**Lemma 4.2:** *For every $\nu \in \mathcal{M}_0$, $X_t$ has a càdlàg version under $P^1_{r,\nu}$.*

**Proof:** It is well-known that (2.11) implies that the total population is finite a.s., which implies that $X$ is càdlàg because $\xi$ is. A rigorous proof of this fact can be based on Neveu's formalism of Galton-Watson trees [Ne86] (cf. [Sc92], Lemma 11). □

**Lemma 4.3:** *For each $f \in b\mathcal{B}^+$, $w(r,x) := P^1_{r,\delta_x}[\langle f, X_t\rangle]$ solves the equation*

$$w(r,x) = \Pi_{r,x}\bigg[f(\xi_t) + \int_r^t w^s(\xi_s)\cdot\alpha^s(\xi_s)\,K(ds)\bigg], \tag{4.2}$$

*where $\alpha(s,x) := \int n\, q(s,x,dn) - 1 = \frac{d}{dz}|_{z=0}\widetilde{\varphi}(s,x,z) \le 0$. In particular,*

$$0 \le P^1_{r,\nu}\big[\langle f, X_t\rangle\big] \le \Pi_{r,\nu}[f(\xi_t)]. \tag{4.3}$$

**Proof:** Substitute $f$ by $\lambda f$ in (2.7) and (2.8), and differentiate with respect to $\lambda$ at $0+$. □

**Lemma 4.4:** *Suppose $U \subset \mathbb{R}_+ \times E$ is open, and $\tau = \inf\{s \ge r \,|\, (s,\xi_s) \in U\}$. Then*

$$P^1_{r,\nu}\Big[\sup_{r\le s\le t} \delta_s \otimes X_s(U) \ge c\Big] \le \frac{1}{c}\Pi_{r,\nu}[\tau \le t].$$

**Proof:** In the critical case, i.e., if $\alpha(r,x)=0$ for all $r$ and $x$, the assertion follows as Equation (30) of [Sc96] by using the historical branching particle system. In the subcritical case, one can compare the process with a critical branching particle systems which is obtained from the original process by adding particles. □

## 5. The tightness proof

In this section we will prove the following.

**Proposition 5.1:** *Suppose $\beta_n \downarrow 0$ as $n \uparrow \infty$. Then, under the assumptions of Theorem 2.4, the distributions of $\beta_n X_t$ under $P^{\beta_n}_{r,\pi_{\mu/\beta_n}}$ are tight on $D(\mathbb{R}_+,\mathcal{M})$.*

Proposition 5.1 will follow by combining Proposition 7.5 with Lemma 5.4 and Lemma 5.6 below.

In the sequel, will now construct a set $\mathbb{F}$ of functions fulfilling the conditions of Proposition 7.5 below. To this end, let $\mathcal{V}$ be a countable base for the topology of $\mathbb{R}_+\times E$ consisting of closed sets. For $V\in\mathcal{V}$, $\tau_{r,V}=\inf\big\{t\ge r \,|\, (t,\xi_t)\in V\big\}$ denotes the first hitting time of $V$ for the process $\xi$ started at time $r$. For $k\in\mathbb{N}$, define

$$f_{V,k}(r,x)=\Pi_{r,x}\big[e^{k(r-\tau_{r,V})}\big]. \tag{5.1}$$

Then it follows by monotone class arguments that $f_{V,k}$ is measurable. For the ease of notations, we will use the abbreviation

$$T^r_t f(x) := \Pi_{r,x}[f(\xi_t)].$$

**Lemma 5.2:** *Fix $k\in\mathbb{N}$ and $V\in\mathcal{V}$, and let $f=f_{V,k}$. Then*

*(i) $g(r,x):=e^{-2kr}f(r,x)$ is an exit rule, i.e., $T^r_t g^t(x)\le g^r$ and $T^r_t g^t(x)\uparrow g^r$ as $t\downarrow r$.*

*(ii) $0\le T^r_t f^t(x)-f^r(x)\le 2k(t-r)$.*

**Proof:** By the Markov property,

$$\begin{aligned} f(r,x) &= \Pi_{r,x}\Big[\exp\big(k(r-\tau_r)\mathrm{I}_{\{\tau_r\,<\,t\}}\big)\,\Pi_{t,\xi_t}\big[e^{k(t-\tau_t)}\big]\Big] \\ &= e^{k(r-t)}\Pi_{r,x}\Big[\exp\big(k(r-\tau_r)\mathrm{I}_{\{\tau_r\,<\,t\}}\big)\,f(t,\xi_t)\Big]. \end{aligned}$$

Thus we get that

$$e^{2k(r-t)}T^r_t f^t \le f^r \le T^r_t f^t.$$

From here both assertions are easily verified. □

For the proof of our next lemma we will need the assumption that $\xi$ is a Hunt process.

**Lemma 5.3:** *Let $f$ be of the form (5.1), for some $k \in \mathbb{N}$ and $V \in \mathcal{V}$, and let $\mu$ be a Borel probability measure on $E$. Then $f$ is $\Pi_{r,\mu}$-quasi continuous in the sense of Definition 7.2 below.*

**Proof:** By Lemma 5.2 there is a $k \in \mathbb{N}$ such that $M_t := e^{-k(r+t)} f(r+t, \xi_{r+t})$ is a supermartingale under $\Pi_{r,\mu}$ with respect to the filtration $\mathcal{G}_t := \mathcal{F}_{[r,r+t]}$, $t \geq 0$. Since $\xi$ is a right process, $M$ is $\Pi_{r,\mu}$-a.s. right continuous; see Theorem 5.3 of [Ku84] or Theorem 5.1 of [Dy82]. By [DM75], VI.3, all the above still holds if the filtered probability space $(\Pi_{r,\mu}, (\mathcal{G}_t))$ is replaced by its completion $(\overline{\Pi}_{r,\mu}, (\overline{\mathcal{G}}_t))$. Moreover, $t \mapsto M_t$ is $\overline{\Pi}_{r,\mu}$-a.s. càdlàg. In view of Proposition 7.4 it thus suffices to show that

$$M_{t-} = e^{-k(r+t)} f(r+t, \xi_{(r+t)-}) \qquad \text{for all } t > 0,\ \overline{\Pi}_{r,\mu}\text{-a.s.} \tag{5.2}$$

To this end, note first that both $t \mapsto M_{t-}$ and $t \mapsto f(r+t, \xi_{(r+t)-})$ are predictable processes. Indeed, for the first process this is clear by left continuity, and the same argument applies to the second one when $f$ is continuous. The general case then follows by monotone class arguments.

To prove (5.2) it thus remains to show that

$$f(\sigma_n, \xi_{\sigma_n}) \longrightarrow f(\sigma, \xi_{\sigma-}) \qquad \overline{\Pi}_{r,\mu}\text{-a.s. on } \{\sigma < \infty\} \tag{5.3}$$

whenever $\sigma_1, \sigma_2, \ldots$ is a sequence of increasing $(\mathcal{G}_t)$-stopping times and $\sigma = \lim_n \sigma_n$; see [DM75], IV.86. But by Lemma 5.2 (ii) and the above we know that $t \mapsto f(r+t, \xi_{r+t})$ is a uniformly integrable $\overline{\Pi}_{r,\mu}$-submartingale, and hence

$$\lim_n f(\sigma_n, \xi_{\sigma_n}) \leq \lim_n \overline{\Pi}_{r,\mu}[f(\sigma, \xi_\sigma) \,|\, \overline{\mathcal{G}}_{\sigma_n}] = f(\sigma, \xi_\sigma) \qquad \overline{\Pi}_{r,\mu}\text{-a.s. on } \{\sigma < \infty\}.$$

Here we have used martingale convergence and the fact that $(\sigma, \xi_\sigma)$ is a.s. $\sigma\Big(\bigcup_n \overline{\mathcal{G}}_{\sigma_n}\Big)$-measurable, because $\xi$ is quasi-left continuous. Repeating the same argument for the supermartingale $e^{-2k(r+t)} f(r+t, \xi_{r+t})$ we arrive at

$$\lim_n f(\sigma_n, \xi_{\sigma_n}) = f(\sigma, \xi_\sigma) \qquad \overline{\Pi}_{r,\mu}\text{-a.s.}$$

Finally we only have to note that $\xi_\sigma = \xi_{\sigma-}$ $\overline{\Pi}_{r,\mu}$-a.s. by quasi-left continuity, and (5.3) is proved. □

**Lemma 5.4:** *The set*

$$\mathbb{F} = \Big\{ h(t,\mu) = \sum_{i=1}^n a_i \exp\langle -f^t_{V_i,k_i}, \mu\rangle \;\Big|\; n \in \mathbb{N},\ V_i \in \mathcal{V},\ k_i \in \mathbb{N},\ a_i \in \mathbb{Q} \Big\}$$

*is uniformly $\big\{P^\beta_{r,\pi_{\mu/\beta}} \,\big|\, 0 < \beta \leq 1\big\}$-quasi continuous in the sense of Definition 7.2 below and separates the points of $\mathbb{R}_+ \times \mathcal{M}$.*

**Proof:** First observe that $f_{V,k} \to \mathrm{I}_V$ pointwise as $k \uparrow \infty$ by right continuity of $\xi$. Thus a $\mathbb{F}$ separates points.

Next suppose $t > r$ and $\varepsilon > 0$ are given. By Lemma 5.3 and Lemma 7.3 (ii) there are compact subsets $A_1 \subset A_2 \subset \cdots \subset [0,t] \times E$ such that each function $f_{V,k}$ is continuous on every $A_n$ and such that $\Pi_{r,\mu}[\tau_n \le t] \le \varepsilon 2^{-n}/2n$, for all $n \in \mathbb{N}$, if $\tau_n$ denotes the first exit time of $(t, \xi_t)$ from $A_n$. Let

$$K = \bigcap_{n=1}^{\infty} \left\{ (s,\gamma) \in [r,t] \times \mathcal{M} \;\middle|\; \delta_s \otimes \gamma(A_n^c) \le \frac{1}{n}, \text{ and } \gamma(E) \le \frac{2}{\langle 1, \mu \rangle \varepsilon} \right\}.$$

Then $K$ is compact by Prohorov's theorem, and

$$\begin{aligned}
P^{\beta}_{r,\pi_{\mu/\beta}} &\Big[ (s, \beta X_s) \notin K \text{ for some } s \le t \Big] \le \\
&\le P^{\beta}_{r,\pi_{\mu/\beta}} \Big[ \sup_{r \le s \le t} \langle 1, \beta X_s \rangle \ge \frac{2}{\langle 1, \mu \rangle \varepsilon} \Big] + \sum_{n=1}^{\infty} P^{\beta}_{r,\pi_{\mu/\beta}} \Big[ \sup_{r \le s \le t} \delta_s \otimes X_s(A_n^c) \ge \frac{1}{\beta n} \Big] \\
&\le \frac{\varepsilon}{2} + \sum_{n=1}^{\infty} n\beta \int \pi_{\mu/\beta}(d\nu)\, \Pi_{r,\nu}[\tau_n \le t] \\
&\le \varepsilon,
\end{aligned}$$

where we have used Lemma 4.4 and the fact that $\langle 1, X_t \rangle$ is a supermartingale under $P^{\beta}_{r,\pi_{\mu/\beta}}$ by Lemma 4.3.

Next we will show that the restriction of each $f \in \mathbb{F}$ to $K$ is continuous. To this end, suppose that $f = f_{V,k}$, for some $V \in \mathcal{V}$ and $k \in \mathbb{N}$, is given, and that $h(s,\gamma) = \langle f^s, \gamma \rangle$. Assume $(s_k, \gamma_k)$ converges to $(s, \gamma)$ in $K$, and let $\delta > 0$ be given. Pick an integer $n$ larger than than $1/\delta$. By the Tietze extension theorem there is a continuous function $\tilde{f}$ from $\mathbb{R}_+ \times E$ to $[0,1]$ such that $\tilde{f}$ coincides with $f$ on $A_n$. Thus

$$\begin{aligned}
\limsup_{k \uparrow \infty} &\big| \langle f^{s_k}, \gamma_k \rangle - \langle f^s, \gamma \rangle \big| \le \\
&\le \limsup_{k \uparrow \infty} \Big( \big| \langle f^{s_k}, \gamma_k \rangle - \langle \tilde{f}^{s_k}, \gamma_k \rangle \big| + \big| \langle \tilde{f}^{s_k}, \gamma_k \rangle - \langle \tilde{f}^s, \gamma \rangle \big| + \big| \langle \tilde{f}^s, \gamma \rangle - \langle f^s, \gamma \rangle \big| \Big) \\
&\le 2 \Big( \delta_s \otimes \gamma(A_n^c) + \limsup_{k \uparrow \infty} \delta_{s_k} \otimes \gamma_k(A_n^c) \Big) \\
&\le 4\delta.
\end{aligned}$$

Since $\delta$ was arbitrary, each function in $\mathbb{F}$ is $\big\{ P^{\beta}_{r,\pi_{\mu/\beta}} \,\big|\, 0 < \beta \le 1 \big\}$-quasi continuous. Now apply Lemma 7.3 (ii) below. $\Box$

**Lemma 5.5:** *Suppose $T > 0$, $0 < \beta \leq 1$, and $h : \mathbb{R}_+ \times \mathcal{M} \to \mathbb{R}$ is of the form $h(t,\mu) = \sum_{i=1}^n a_i \exp\langle -f_i^t, \mu\rangle$ with $a_i \in \mathbb{R}$ and $f_i : \mathbb{R}_+ \times E \to \mathbb{R}_+$ bounded and measurable. Then there exists a constant $C$ depending only on $h$ and $T$ such that*

$$P^\beta_{r,\nu}\Big[\Big(h(t,\beta X_t) - h(r,\beta\nu)\Big)^2\Big] \leq C\Big(\sum_{i=1}^n \|T_t^r f_i^t - f_i^r\| + k(r,t)\Big)\langle 1, \beta\nu\rangle,$$

*for all $0 \leq r < t \leq T$, where $k$ was introduced in Lemma 4.1.*

**Proof:** During the proof, $c$ will denote a varying generic constant depending only on $h$. Since

$$P^\beta_{r,\nu}\Big[\big(h(t,\beta X_t) - h(r,\beta\nu)\big)^2\Big] \leq c\sum_{i=1}^n P^\beta_{r,\nu}\Big[\Big(e^{-\langle f_i^t, \beta X_t\rangle} - e^{-\langle f_i^r, \beta\nu\rangle}\Big)^2\Big],$$

we may assume that $n = 1$ and $a_1 = 1$, and we may suppress the index $n = 1$ in the sequel. Then let $v_\beta(r,x) = P^\beta_{r,\delta_x}\big[\exp\langle -f^t, \beta X_t\rangle\big]$, and define $w_\beta$ as $v_\beta$ but with $f$ replaced by $2f$. We have

(5.4)

$$\begin{aligned} P^\beta_{r,\nu}\Big[\big(h(t,\beta X_t) - h(r,\beta\nu)\big)^2\Big] &= \\ &= e^{\langle \log w^r_\beta, \nu\rangle} - 2e^{\langle \log v^r_\beta - \beta f^r, \nu\rangle} + e^{-\langle 2\beta f^r, \nu\rangle} \\ &\leq \Big| e^{\langle \log w^r_\beta, \nu\rangle} - e^{\langle \log v^r_\beta - \beta f^r, \nu\rangle}\Big| + \Big| e^{\langle \log v^r_\beta, \nu\rangle} - e^{-\langle \beta f^r, \nu\rangle}\Big| \\ &\leq \big\langle \big| \log w^r_\beta - \log v^r_\beta - \log e^{-\beta f^r}\big|, \nu\big\rangle + \big\langle \big| \log v^r_\beta - \log e^{-\beta f^r}\big|, \nu\big\rangle \\ &\leq c\Big(\big\langle |w^r_\beta - e^{-2\beta f^r}|, \nu\big\rangle + 2\big\langle |v^r_\beta - e^{-\beta f^r}|, \nu\big\rangle\Big). \end{aligned}$$

Now

$$\|v^r_\beta - e^{-\beta f^r}\| \leq \|v^r_\beta - T^r_t e^{-\beta f^t}\| + \|T^r_t e^{-\beta f^t} - e^{-\beta f^r}\|. \tag{5.5}$$

Next we will estimate separately both terms on the right hand side of (5.5). The second term can be estimated by

$$\|T^r_t e^{-\beta f^t} - e^{-\beta f^r}\| \leq c\beta \|T^r_t f^t - f^r\|. \tag{5.6}$$

As for the first term note that by Lemma 4.3

$$v_\beta(r,x) = 1 - P^\beta_{r,\delta_x}\big[1 - e^{-\beta\langle f^t, X_t\rangle}\big] \geq 1 - \beta\cdot P^\beta_{r,\delta_x}\big[\langle f^t, X_t\rangle\big] \geq 1 - \beta\cdot\|f^t\|. \tag{5.7}$$

Since by assumption $\psi_\beta(s,y,z) = \widehat{\phi}_\beta(s,y,1-\beta z)/\beta^2$ converges locally uniformly as $\beta \downarrow 0$, we infer from (5.7) that there is a constant $L$ such that, for all $\beta \in (0,1]$,

$$0 \leq \frac{1}{\beta}\widehat{\phi}_\beta(s,y,v^s_\beta(y)) \leq \beta\cdot L \qquad \text{uniformly in } y \text{ and } s \leq T.$$

Thus

$$0 \leq v_\beta(r,x) - T^r_t e^{-\beta f^t}(x) = P^\beta_{r,\delta_x}\Big[\frac{1}{\beta}\int_r^t \widehat{\phi}_\beta(s,\xi_s,v^s_\beta(\xi_s))\,K(ds)\Big] \leq \beta L\cdot\Pi_{r,x}\big[K[r,t]\big].$$

Together with (5.6), (5.5), and (5.4) the assertion follows. □

**Lemma 5.6:** *Suppose $\beta_n \downarrow 0$. Then, for any $h \in \mathbb{F}$, the distributions of the process $t \mapsto h(t, \beta_n X_t)$ under $P^{\beta_n}_{r,\pi_{\mu/\beta_n}}$, $n \in \mathbb{N}$, are tight on $D(\mathbb{R}_+; \mathbb{R})$.*

**Proof:** Fix $r$, $T > 0$, and let for the moment $\mathcal{G}_t$ denote the $\sigma$-field generated by $X_s$, $r \le s \le t$. Then by Lemma 5.5 and Lemma 5.2 (ii) there is a constant $C$ such that, for all $u \in [0, \delta]$ and $t \in [r, T]$,

$$P^{\beta}_{r,\pi_{\mu/\beta}}\Big[\Big(h(t+u, \beta X_{t+u}) - h(t, \beta X_t)\Big)^2 \ \Big|\ \mathcal{G}_t\Big]$$
$$= P^{\beta}_{r,X_t}\Big[\Big(h(t+u, \beta X_{t+u}) - h(t, \beta X_t)\Big)^2\ \Big]$$
$$\le C(\delta + k(t, t+u))\langle 1, \beta X_t\rangle.$$

Now let $\tilde{k}(\delta) := \sup_{0 \le t \le T} k(t, t+\delta)$ and let $M^\beta - A^\beta$ be the Doob-Meyer decomposition of the supermartingale $\langle 1, X_t\rangle$ under $P^{\beta}_{r,\pi_{\mu/\beta}}$. Then (7.7) holds in our setting if $\gamma_\beta(\delta)$ is defined by

$$\gamma_\beta(\delta) = C(\delta + \tilde{k}(\delta))\beta M^\beta_T.$$

Moreover, for all $\beta \in (0, 1]$,

$$P^{\beta}_{r,\pi_{\mu/\beta}}\big[\gamma_\beta(\delta)\big] = C(\delta + \tilde{k}(\delta))\langle 1, \mu\rangle \longrightarrow 0 \qquad \text{as } \delta \downarrow 0$$

by Lemma 4.1. Thus the assertion follows from Proposition 7.6. □

## 6. Proof of the regularity results.

**Proof of Theorem 2.4:** Under the assumptions of the theorem the functions $v_\beta$ defined in (2.9), (2.10) converge uniformly in $(r, x) \in [0, T] \times E$ and in $f \in b\mathcal{B}^+$ with $\|f\| \le c$ to the unique positive solution of (2.2). The proof is the same as in [Dy94], pp. 56-58. Using the Markov property and induction it follows that the finite-dimensional distributions of the approximating branching particle systems converge weakly. Hence Theorem 2.4 follows from Proposition 5.1. □

**Proof of Theorem 2.3:** It suffices to prove the existence of a càdlàg version up to time $T > 0$. To this end, let $h$, $\xi^h$, $\psi_h$, $K^h$, and $\mu^h$ be as in the proof of Theorem 2.2. According to [Dy91a], for every superprocess satisfying the assumptions of Theorem 2.1, there exists a sequence of approximating branching particle systems satisfying the assumptions of Theorem 2.4. Therefore we can assume that the $(\xi^h, \psi_h, K^h)$-superprocess is already càdlàg. Observe that $h$ is $\Pi^h_{r,\mu^h}$-quasi continuous (cf. the proof of Lemma 5.3 and note that the proof is even easier in this case). By Lemma 7.3 and Lemma 5.4 and its proof one easily deduces that $\Phi(t, \mu) := h^{-1}(t, x)\,\mu(dx)$ is a $P^h_{r,\mu^h}$-quasi homeomorphism from $\mathcal{M}$ to $\mathcal{M}^\varrho$ (with the obvious definition of a quasi-homeomorphism). Thus the process $h^{-1}(t,x)\,X_t(dx)$ is càdlàg $P^h_{r,\mu^h}$-a.s. in the same way, the right property follows from the right property in the case $\varrho \equiv 1$, where it is well-known (see e.g. [Ku94]). $\square$

**Proof of Corollary 2.5:** we leave it to the reader to check that the assumptions of Theorem 2.4 are also fulfilled in case of the historical branching particle systems. As for (ii), $f(t,\widehat{x}) := g(\pi_t(\widehat{x}))$ has the property that $t \mapsto f(t,\widehat{\xi}_t)$ is càdlàg and that $t \mapsto f(t,\widehat{\xi}_{t-})$ is càglàd $\widehat{\Pi}_{r,\hat{\mu}}$-a.s., for all continuous functions $g$ on $E$ and all $\widehat{\mu}$. Thus $f$ is $\widehat{\Pi}_{r,\hat{\mu}}$-quasi continuous, and it follows that $(t,\widehat{\nu}) \mapsto \widehat{\nu}\circ\pi_t^{-1}$ is a $\widehat{P}_{r,\hat{\mu}}$-quasi continuous map. $\square$

## 7. Appendix: Skorohod space and quasi-continuity

In this section, $S$ is a separable and metrizable topological space, $I$ equals either $\mathbb{R}_+$ or $[0,T]$ with $T>0$ finite, and $(Z_t)_{t\in I}$ denotes the coordinate process on the Skorohod space $D(I;S)$, i.e. $Z_t(\omega)=\omega(t)$ for all $\omega\in D(I;S)$ and $t\in I$. As usual, $D(I;S)$ will be endowed with the Skorohod topology and the corresponding Borel $\sigma$-algebra, which is also generated by $(Z_t)_{t\in I}$. See [EK86] or [Ja86]. By saying a subset of a metric space is analytic we mean that it is analytic with respect to the paving formed by all Borel sets; see [DM75]. A co-analytic set is the complement of an analytic set.

**Proposition 7.1:**

*(i) Suppose $C$ is a co-analytic subset of $S$. Then the set $D(I;C)$ (defined with respect to the relative topology on $C\subset S$) is a co-analytic subset of $D(I;S)$.*

*(ii) If $S$ is a metrizable co-Souslin space, then also $D(I;S)$ is a metrizable co-Souslin space.*

**Proof:** To prove (i), note first that $(t,w)\mapsto w(t)$ and $(t,w)\mapsto w(t-)$ are both measurable mappings from $I\times D(I;S)$ to $S$. By [DM75], III.11, we hence get that

$$B := \Big\{(t,w)\;\Big|\; w(t)\in C^C \;\text{ or } w(t-)\in C^C\Big\}$$

is an analytic subset of $I\times D(I;S)$. Therefore also the projection $\overline{B}$ of $B$ onto $D(I;S)$ is analytic in $D(I;S)$ by [DM75], III.13.3. But the complement $\overline{B}^C$ of $\overline{B}$ in $D(I;S)$ is just the set of all paths from $I$ into $C$ which are right continuous and have left limits in $C$. I.e., $\overline{B}^C = D(I;C)$.

In order to prove (ii) we can assume that there is compact metric space $K$ such that $S$ is a co-analytic subset of $K$ and that $S$ carries the topology relative to this embedding. Then, by (i), $D(I,S)$ is a co-analytic subset of the Polish space $D(I;K)$. But a standard characterization of the Skorohod topology (see, for instance, Proposition 6.5 in Chapter 3 of [EK86]) shows that the topology induced on $D(I;S)$ by the inclusion $D(I;S) \subset D(I;K)$ coincides with the Skorohod topology. Thus $D(I;S)$ is a metrizable co-Souslin space. □

**Definition 7.2:** *Suppose $\mathcal{P}$ is a set of Borel probability measures on $D(I,S)$, and $\mathbb{F}$ is a family of functions on $I \times S$. We will say that* $\mathbb{F}$ is uniformly $\mathcal{P}$-quasi continuous, *if there is an increasing sequence $(K_n)_{n\in\mathbb{N}}$ of compact subsets of $I \times S$ such that*

$$\sup_{P\in\mathcal{P}} P\Big[(s, Z_s) \notin K_n \text{ for some } s \in [0,n] \cap I\Big] \longrightarrow 0 \qquad \text{as } n \uparrow \infty. \tag{7.1}$$

$$\text{Each function } f \in \mathbb{F} \text{ is continuous on every } K_n. \tag{7.2}$$

*If both $\mathbb{F}$ and $\mathcal{P}$ consist only of single elements $f$ and $P$, respectively, we will just say that $f$ is $P$-quasi continuous instead of $\{f\}$ is uniformly $\{P\}$-quasi continuous.*

The following lemma is easy to prove.

**Lemma 7.3:**

*(i) If $\mathbb{F}$ is uniformly $\mathcal{P}$-quasi continuous, then $\mathbb{F}$ is also uniformly $\overline{\mathcal{P}}$-quasi continuous, where $\overline{\mathcal{P}}$ is the weak closure of $\mathcal{P}$.*

*(ii) Suppose $\mathcal{P}$ is a family of Borel probability measures on $D(I;S)$, and $\mathbb{F}$ is a countable set of functions on $I \times S$. If $\{f\}$ is uniformly $\mathcal{P}$-quasi continuous for every $f \in \mathbb{F}$, then $\mathbb{F}$ is uniformly $\mathcal{P}$-quasi continuous.*

We have the following criterion for a function $f$ to be $P$-quasi continuous. It was proved first by Le Jan [Le83] for a locally compact space $S$.

**Proposition 7.4:** *Suppose that $S$ is a metrizable co-Souslin space, $f$ is a measurable function on $I \times S$, and $P$ is a Borel probability measure on $D(I;S)$. Then the following two conditions are equivalent.*

*(i) $f$ is $P$-quasi continuous.*

*(ii) The processes $t \mapsto f(t, Z_t)$ and $t \mapsto f(t, Z_{t-})$ are $P$-a.s. càdlàg and càglàd respectively.*

**Proof:** Suppose condition (ii) holds. We may assume that $f$ is bounded, and we write $\widetilde{S}$ for $I \times S$, $\widetilde{Z}$ for the space-time process $\widetilde{Z}_t = (t, Z_t)$, $t \in I$, and $\widetilde{P}$ for the image measure $P \circ \widetilde{Z}^{-1}$. Following [Le83], we introduce a norm $\|g\|_P$ for a measurable function $g$ on $\widetilde{S}$ as follows.

$$\|g\|_P = P\Big[\sup_{t\in I} 2^{-t}|g(\widetilde{Z}_t)| \vee |g(\widetilde{Z}_{t-})|\Big].$$

One can show as in [MR92], pp. 134, 135 that there is a sequence $(f_n)$ of bounded and continuous functions on $I \times S$ such that $\|f - f_n\|_P \to 0$. We may assume that $\|f_{n+1} - f_n\|_P \le 2^{-3n}$, and that the $f_n$ are uniformly bounded. Writing $I_n$ for $I \cap [0, n]$ we let

$$F_n = \Big\{(t,x) \in I_n \times S \;\Big|\; |f_{n+1}(t,x) - f_n(t,x)| > 2^{-n}\Big\} \quad \text{and} \quad B_N = \bigcup_{n \ge N} F_n.$$

Then, due to Proposition 7.1 (i) and [DM75], III.33, the spaces $D(I_n; F_n^C)$ and $D(I_N; B_N^C)$ are co-analytic and hence $\widetilde{P}$-measurable subsets of $D(I; \widetilde{S})$. By construction we get that

$$\begin{aligned}\widetilde{P}\Big[D(I_N; B_N^C)^C\Big] &\le \sum_{n \ge N} P\Big[\exists\, t \le n \text{ such that } (t, Z_t) \in F_n \text{ or } (t, Z_{t-}) \in F_n\Big] \\ &\le \sum_{n \ge N} 2^{2n} \|f_{n+1} - f_n\|_P \le 2^{1-N}.\end{aligned} \tag{7.3}$$

According to Proposition 7.1 (ii) and [DM75], III.38, $\widetilde{P}$ is an inner regular measure on $D(I; \widetilde{S})$. Therefore there exist compact sets $\mathcal{K}_N \subset D(I_N; B_N^C)$ such that $\widetilde{P}[\mathcal{K}_N] \ge 1 - 2^{2-N}$. But this implies that there are compact sets $K_N \subset B_N$ such that

$$P\Big[\exists\, t \le N \text{ such that } (t, Z_t) \notin K_N\Big] \le 2^{2-N} \longrightarrow 0 \qquad (N \uparrow \infty) \tag{7.4}$$

(cf. Proposition 1.6 (iv) in [Ja86]). Note that $(f_n)$ converges uniformly on each $K_N$. Therefore the measurable function $\widetilde{f}$ defined by

$$\widetilde{f}(t,x) = \begin{cases} \lim_n f_n(t,x) & \text{if } (t,x) \in \bigcup_N K_N, \\ 0 & \text{otherwise,} \end{cases}$$

is $P$-quasi continuous. Moreover, $\|\widetilde{f} - f_n\|_P \to 0$, since the $f_n$ are uniformly bounded. Hence we conclude that $\|f - \widetilde{f}\|_P = 0$. Consider the set $A = \big\{(t,x) \in \overline{S} \mid f(t,x) \neq \widetilde{f}(t,x)\big\}$. By a similar line of reasoning as above we get the existence of a sequence $(\widetilde{K}_N)$ of compact subsets of $A$ such that (7.4) holds with $K_N$ replaced by $\widetilde{K}_N$, and hence also for $K_N \cap \widetilde{K}_N$. But $f$ is continuous on every set $K_N \cap \widetilde{K}_N$, i.e. $f$ is $P$-quasi continuous. Thus we have proved (ii) $\Rightarrow$ (i). The converse statement is trivial. $\square$

Any continuous function $f$ on $I \times S$ induces a mapping $f_*$ as follows.

$$\begin{aligned} f_* : D(I;S) &\longrightarrow D(I;\mathbb{R}) \\ \omega &\longmapsto \Big(t \mapsto f(t,\omega(t))\Big). \end{aligned} \tag{7.5}$$

If the function $f$ is not continuous, but $P$-quasi continuous for some $P$ on $D(I;S)$, then $f_*$ in the sense of (7.5) is still defined $P$-a.s. This allows us to generalize the well known Jakubowski criterion for tightness on $D(I;S)$ in the following way.

**Proposition 7.5:** *Suppose $\mathcal{P}$ is a family of Borel probabilities on $D(I;S)$. Then the following are equivalent.*

1. *$\mathcal{P}$ is tight.*
2. *There exists a family $\mathbb{F}$ of $\mathcal{P}$-quasi continuous functions closed under addition and separating the points of $I \times S$ such that $(P \circ f_*^{-1})_{P \in \mathcal{P}}$ is tight on $D(I;\mathbb{R})$, for each $f \in \mathbb{F}$.*

The proof of Proposition 7.5 is essentially the same as the one of Theorem 4.6 in [Ja86], and therefore it is omitted. We close this section with recalling the following tightness criterion on $D(I;\mathbb{R})$, which is a simplified version of Theorem 3.8.6 of [EK86]. Recall the notation $I_k = I \cap [0,k]$.

**Proposition 7.6:** *Suppose $P_1, P_2, \ldots$ are Borel probability measures on $D(I;\mathbb{R})$ such that for all $\varepsilon$, $k > 0$ there is a $L > 0$ such that*

$$\sup_n P_n\Big[\sup_{t \in I_k} |Z_t| > L\Big] \le \varepsilon, \tag{7.6}$$

*and such that, for each $k$, there are random variables $\gamma_n(\delta)$, $0 < \delta < 1$, on $D(I;\mathbb{R})$ with*

$$\lim_{\delta \downarrow 0} \limsup_n P_n[\,\gamma_n(\delta)\,] = 0$$

*and with*

$$P_n\Big[(Z_{t+u} - Z_t)^2 \;\Big|\; Z_s,\, s \le t\Big] \le P_n\Big[\,\gamma_n(\delta) \;\Big|\; Z_s,\, s \le t\Big] \tag{7.7}$$

*for all $n \in \mathbb{N}$, $t \in I_k$, and $0 \le u \le \delta$ such that $t + u \in I_k$. Then $\{P_1, P_2, \ldots\}$ is tight.*

**References:**

[Da93] Dawson, D A. Measure-valued Markov processes. École d'Été de Probabilités de Saint Flour, 1991, Lecture Notes in Math. 1541, Springer, 1993.

[DF97] Dawson, D. A. and Fleischmann, K. A continuous super-Brownian motion in a super-Brownian medium. J. Theoret. Probab. 10 (1997), no. 1, 213–276.

[DFL97] Dawson, D. A., Fleischmann, K. and Leduc, G. Continuous dependence of a class of superprocesses on branching parameters, and applications. WIAS preprint No. 317 (1997).

[DP91] Dawson, D A. and Perkins, E. Historical processes. Mem. Amer. Math. Soc. 454, 1991.

[DM75] Dellacherie, C. and Meyer, P.-A. Probabilités et Potentiel. Hermann, Paris, 1975.

[Dy82] Dynkin, E. B. Regular Markov processes. In: Dynkin, E. B. Markov processes and related problems in analysis. Cambridge University Press, London, 1982.

[Dy91a] Dynkin, E. B. Branching particle systems and superprocesses. Ann. Probab. 19 (1991), 1157–1194.

[Dy91b] Dynkin, E. B. Path processes and historical superprocesses. Probab. Theory Related Fields 90 (1991), 1–36.

[Dy94] Dynkin, E. B. An introduction to branching measure-valued processes. CRM Monograph Series, 6. American Mathematical Society, Providence, 1994.

[EK86] Ethier, S. N. and Kurtz, T. G. Markov processes — characterization and convergence. Wiley, New York, 1986.

[Fi88] Fitzsimmons, P. J. Construction and regularity of measure-valued Markov branching processes. Israel J. Math. 64 (1988), 337–361.

[FM97] Fleischmann, K. and Mueller, C. A super-Brownian motion with a locally infinite catalytic mass. Probab. Theory Related Fields 107 (1997), 325–357.

[Is86] Iscoe, I. A weighted occupation time for a class of measure-valued critical branching processes. Probab. Theory Relat. Fields 71 (1986), 85–116.

[Ja86] Jakubowski, A. On the Skorokhod topology. Ann. Inst. H. Poincaré Probab. Statist. 22 (1986), no. 3, 263–285.

[Ko99] Kolesnikov, A. V. Über die topologischen Eigenschaften des Skorohod-Raumes. To appear: Theory Probab. Appl. (1999).

[Ku84] Kuznetsov, S. E. Nonhomogeneous Markov processes. J. Soviet Mat. 25 (1984), 1380-1498.

[Ku94] Kuznetsov, S. E. Regularity properties of a supercritical superprocess. The Dynkin Festschrift, 221–235, Progr. Probab., 34, Birkhäuser, Boston, 1994.

[Le83] Le Jan, Y. Quasi-continuous functions and Hunt processes. J. Math. Soc.

Japan 35 (1983), 37-42.

[MR92] Ma, Z. M. and Röckner, M. Introduction to the theory of (nonsymmetric) Dirichlet forms. Springer-Verlag, Berlin, 1992.

[Ne86] Neveu, J. Arbres et processus de Galton–Watson. Ann. Inst. H. Poincaré Probab. Statist. 22 (1986), 199–207.

[R-C86] Roelly-Coppoletta, S. A criterion of convergence of measure-valued processes: application to measure branching processes. Stochastics 17 (1986), 43-65.

[Sc92] Schied, A. Zur Konstruktion maßwertiger Verzweigungsprozesse. Diplomarbeit, Univ. Bonn (1992).

[Sc96] Schied, A. Sample path large deviations for super-Brownian motion. Probab. Theory Related Fields 104 (1996), 319–347.

[Wa68] Watanabe, S. A limit theorem of branching processes and continuous state branching processes. J. Math. Kyoto Univ. 8 (1968), 141–167.

Alexander Schied, Institut für Mathematik, Humboldt-Universität, Unter den Linden 6, D-10099 Berlin, Germany

schied@mathematik.hu-berlin.de